# $L^p$-theory for fractional gradient PDE with VMO coefficients

Armin Schikorra[*1], Tien-Tsan Shieh[†2], and Daniel Spector[‡3]

[1]Department of Mathematics, University of Basel, Basel, Switzerland
[2]Department of Mathematics, National Central University, Chung Li, Taiwan
[3]Department of Applied Mathematics, National Chiao Tung University, Hsinchu, Taiwan

**Abstract**

In this paper, we prove $L^p$ estimates for the fractional derivatives of solutions to elliptic fractional partial differential equations whose coefficients are $VMO$. In particular, our work extends the optimal regularity known in the second order elliptic setting to a spectrum of fractional order elliptic equations.

## 1 Introduction

In his 1959 paper on some composition formulas for vector-valued potentials, J. Horváth introduced [7, p. 434] the differential object

$$D^s u := D I_{1-s} u. \tag{1.1}$$

Here, $s \in (0,1)$ and $I_{1-s}$ is the Riesz potential of order $1-s$.

This object was subsequently termed the *Riesz fractional gradient* by the second and third author in [13], where it was utilized to generalize divergence form elliptic partial differential equations from the second order setting to that of differential order $2s \in (0,2)$. In particular, assuming that $A$ is uniformly elliptic, i.e.

$$\lambda |\xi|^2 \leq A(x)\xi \cdot \xi \leq \Lambda |\xi|^2, \tag{1.2}$$

[*]armin.schikorra@unibas.ch, A.S. supported by SNF
[†]tshieh@math.ncu.edu.tw, T.-T. S. supported by MOST 103-281-M-008-066
[‡]dspector@math.nctu.edu.tw, D.S. supported by MOST 103-2115-M-009-016-MY2

for all $x, \xi \in \mathbb{R}^N$ and some $0 < \lambda \leq \Lambda < +\infty$, the authors showed that given $\varphi \in H^s(\mathbb{R}^N)$ and $g \in L^2(\Omega)$ there exists $u \in H^s(\mathbb{R}^N)$ that satisfies

$$\int_{\mathbb{R}^N} A(x) D^s u(x) \cdot D^s v(x) \, dx = \int_{\mathbb{R}^N} gv \tag{1.3}$$

for all $v \in C_c^\infty(\mathbb{R}^N)$ and $u = \varphi$ in $\mathbb{R}^N \setminus \Omega$. Here, $\Omega \subset \mathbb{R}^N$ is open and bounded, $N \geq 2$, and

$$H^s(\mathbb{R}^N) := \{u \in L^2(\mathbb{R}^N) : D^s u \in L^2(\mathbb{R}^N; \mathbb{R}^N)\},$$

which coincides with any standard definition of the fractional Sobolev space.

One observes that when $s = 1$ and the boundary of $\Omega$ is sufficiently nice, the equation (1.3) agrees with the weak formulation of a divergence form elliptic PDE, since prescribing $u$ on the complement gives rise to a trace that would be a more standard way to frame the existence. Meanwhile for $s \in (0, 1)$ one obtains a family of fractional partial differential equations with analogous structure. The interest in generalizing partial differential equations via (1.1) is two-fold. Firstly, that one should be concerned with non-integer order differential objects can be simply explained by quoting Sobolev and Nikol'skiĭ's 1963 paper (who even implicitly consider (1.1), see [12, p. 148]) where they note that "an imbedding theory containing only derivatives of integral order is incomplete and imperfect." Secondly, the structure of (1.1) closely resembles the gradient and therefore such a generalization preserves the structural properties of the equation, a point which we will return to later. This aspect has been important in the development of $L^1$ fractional Sobolev inequalities in terms of (1.1) in [11], as such inequalities are known to be false for the fractional Laplacian.

In this paper we continue to develop this perspective of classical equations as a part of a continuous spectrum. In particular, we take the first step in addressing for this class of equations a question of fundamental importance in the second order case, that of regularity. As there are a number of possible assumptions one can make to investigate the question of regularity of $u$ that satisfies (1.3), let us further describe the hypothesis of interest to us. In addition to the ellipticity condition (1.2), we will assume $A$ is of vanishing mean oscillation.

**Definition 1.1** *We define the semi-norm (on the space of functions of bounded mean oscillation)*

$$[\varphi]_{BMO} := \sup_Q \fint_Q |\varphi - \fint_Q \varphi|.$$

*Then we define the space of functions of* vanishing mean oscillation *by*

$$VMO(\mathbb{R}^N) := \overline{\{C_c^\infty(\mathbb{R}^N)\}}^{[\cdot]_{BMO}}.$$

The main result of this paper is the following theorem on the regularity of such equations with VMO coefficients.

**Theorem 1.2** *Suppose that $A \in VMO(\mathbb{R}^N; \mathbb{R}^{N \times N})$ satisfies (1.2), that $G \in L^p(\mathbb{R}^N; \mathbb{R}^N)$ for some $1 < p < +\infty$ and $u \in H^s(\mathbb{R}^N)$ satisfies*

$$\int_{\mathbb{R}^N} A(x) D^s u(x) \cdot D^s v(x) \, dx = \int_{\mathbb{R}^N} G \cdot D^s v \tag{1.4}$$

*for all* $v \in C_c^\infty(\Omega)$. *Then* $D^s u \in L^p_{loc}(\Omega)$ *and for any* $K \subset\subset \Omega$ *there exists a constant* $C = C(K, \Omega, A, s, p) > 0$ *such that*

$$\|D^s u\|_{L^p(K;\mathbb{R}^N)} \leq C \left( \|G\|_{L^p(\mathbb{R}^N;\mathbb{R}^N)} + \|(-\Delta)^{\frac{s}{2}} u\|_{L^2(\mathbb{R}^N)} \right).$$

Here, $(-\Delta)^{\frac{s}{2}} u$ denotes the fractional Laplacian of $u$ of order $s$, which can be defined as a Fourier multiplier with symbol $(2\pi|\xi|)^s$, see [14, p. 117]. The fractional Laplacian is related to the fractional gradient via the identity

$$D^s u \equiv R(-\Delta)^{\frac{s}{2}} u, \tag{1.5}$$

for $s \in (0, 1)$ and $u$ with sufficient smoothness and integrability, and where $R = DI_1$ is the vector-valued Riesz transform. In the sequel we take (1.5) as our definition of $D^s u$, which enables us to include the classical case $s = 1$ (and more generally $s > 1$ though one loses the interpretation of a fractional gradient in this range).

Our proof is based on the beautiful technique of Iwaniec and Sbordone, introduced in [8] for $u$ satisfying (1.4) with $v \in C_c^\infty(\mathbb{R}^N)$ and $s = 1$. We recall that in this setting they had shown [8, p. 186] that (1.4) has exactly one (up to a constant) solution with the estimate

$$\|Du\|_{L^p(\mathbb{R}^N;\mathbb{R}^N)} \leq C \|G\|_{L^p(\mathbb{R}^N;\mathbb{R}^N)}.$$

Comparing this with our result, one sees that the preservation of structure in the equation results in regularity that is completely analogous to the well-studied elliptic theory.

As a consequence of this result we can return to the question of regularity of solutions to (1.3). In particular, one can transform equation (1.3) into (1.4) by defining $G = I_s R g$ (where one extends $g$ by zero outside $\Omega$), since one has

$$\int_{\mathbb{R}^N} g v \, dx = \int_{\mathbb{R}^N} I_s R g \cdot R(-\Delta)^{\frac{s}{2}} v \, dx$$
$$= \int_{\mathbb{R}^N} G \cdot D^s v \, dx$$

for $v \in C_c^\infty(\mathbb{R}^N)$ and $g \in L^2(\mathbb{R}^N)$. The assumption $g \in L^2(\Omega)$ then implies that $G \in L^{2N/(N-2s)}(\mathbb{R}^N; \mathbb{R}^N)$, and so our result allows us to conclude that for the solution to (1.3) we have for every $K \subset\subset \Omega$ the estimate

$$\|D^s u\|_{L^{2N/(N-2s)}(K;\mathbb{R}^N)} \leq C \left( \|g\|_{L^2(\Omega)} + \|(-\Delta)^{\frac{s}{2}} u\|_{L^2(\mathbb{R}^N)} \right).$$

When $s = 1$ this localizes the result of Iwaniec and Sbordone and can be compared with a result of Di Fazio in [5] (who in fact obtains regularity up to the boundary).

## 2 Estimates and proof of the Main Result

The main tool we utilize is the following result of Iwaniec and Sbordone [8, see p. 187, 201-206].

**Theorem 2.1 (Iwaniec, Sbordone)** *Let $A \in VMO \cap L^\infty(\mathbb{R}^N; \mathbb{R}^{N \times N})$ satisfy (1.2). Then for all $1 < q < +\infty$, the operator*

$$T := R_i A_{ij} R_j : L^q(\mathbb{R}^N) \to L^q(\mathbb{R}^N)$$

*is invertible, and moreover, there exists $C = C(A, q) > 0$ such that*

$$\|f\|_{L^q(\mathbb{R}^N)} \leq C \|Tf\|_{L^q(\mathbb{R}^N)} \tag{2.1}$$

*for all $f \in L^q(\mathbb{R}^N)$.*

From this we obtain the localization:

**Proposition 2.2** *Let $A, T$ as in Theorem 2.1. Then for any $\Omega_1, \Omega_2$ open and bounded with $\Omega_1 \subset\subset \Omega_2$, $2 < q < +\infty$, there exists $C = C(A, q, \Omega_1, \Omega_2) > 0$ such that*

$$\|f\|_{L^q(\Omega_1)} \leq C(\|Tf\|_{L^q(\Omega_2)} + \|f\|_{L^2(\mathbb{R}^N)})$$

*for all $f \in L^2(\mathbb{R}^N)$.*

Before proving Proposition 2.2, let us recall the following commutator estimate, whose proof we provide for the convenience of the reader.

**Proposition 2.3** *Let $b, f : \mathbb{R}^N \to \mathbb{R}$ and define the commutator $\mathcal{C}(b, R_i)[f]$ by*

$$\mathcal{C}(b, R_i)[f] := bR_i[f] - R_i[bf],$$

*where $R_i$ is the $i$-th Riesz transform. If $b$ is Lipschitz, then*

$$\|\mathcal{C}(b, R_i)[f]\|_{L^p(\mathbb{R}^N)} \leq C[b]_{\mathrm{Lip}(\mathbb{R}^N)} \|I_1|f|\|_{L^p(\mathbb{R}^N)}.$$

**Proof.** Since

$$R_i g(x) = c_N \int_{\mathbb{R}^N} \frac{x_i - z_i}{|x - z|^{N+1}} g(z) \, dz,$$

we have

$$\mathcal{C}(b, R_i)[f](x) = c_N \int_{\mathbb{R}^N} \frac{x_i - z_i}{|x - z|^{N+1}} (b(x) - b(z)) f(z) \, dz,$$

and consequently,

$$|\mathcal{C}(b, R_i)[f](x)| \leq c_N [b]_{\mathrm{Lip}(\mathbb{R}^N)} \int_{\mathbb{R}^N} |x - z|^{-N+1} |f|(z) \, dz = C[b]_{\mathrm{Lip}(\mathbb{R}^N)} I_1|f|(x).$$

∎

**Proof of Proposition 2.2.** Let $\eta \in C_0^\infty(\Omega_2)$ be a usual cutoff function, i.e. $\eta \geq 0$ and $\eta \equiv 1$ on a neighbourhood of $\Omega_1$. From (2.1) we have

$$\|f\|_{L^q(\Omega_1)} \leq \|\eta f\|_{L^q(\mathbb{R}^N)} \leq C \|T(\eta f)\|_{L^q(\mathbb{R}^N)}.$$

Let us now recall the definition of the commutator of an operator $T$ and two functions $b, f$ (which can be thought of as the error term to a product rule). We have

$$\mathcal{C}(b, T)[f] := bT[f] - T[bf].$$

Then we continue the preceding estimate as follows. For $\operatorname{supp}\eta \subset\subset K_0 \subset\subset L_1 \subset\subset \Omega_2$ and denoting $\chi_{L_1}$ the characteristic function of $L_1$, we estimate

$$\begin{aligned}
\|T(\eta f)\|_{L^q(\mathbb{R}^N)} =& \|T(\eta \chi_{L_1} f)\|_{L^q(\mathbb{R}^N)}\\
\leq & \|\eta T(\chi_{L_1} f)\|_{L^q(\mathbb{R}^N)} + \|\mathcal{C}(\eta, T)[\chi_{L_1} f]\|_{L^q(\mathbb{R}^N)}\\
\leq & \|T(\chi_{L_1} f)\|_{L^q(K_0)} + \|\mathcal{C}(\eta, T)[\chi_{L_1} f]\|_{L^q(\mathbb{R}^N)}\\
\leq & \|T(f)\|_{L^q(K_0)} + \|T(\chi_{L_1^c} f)\|_{L^q(K_0)} + \|\mathcal{C}(\eta, T)[\chi_{L_1} f]\|_{L^q(\mathbb{R}^N)}\\
=: & \|T(f)\|_{L^q(K_0)} + I + II.
\end{aligned}$$

Note that in the above with our $T$ we have

$$\mathcal{C}(\eta, T)[\chi_{L_1} f] = R_i A_{ij} [\mathcal{C}(\eta, R_j)[\chi_{L_1} f]] + \mathcal{C}(\eta, R_i)[A_{ij} R_j(\chi_{L_1} f)].$$

As for $I$, since the supports of $L_1^c$ and $K_0$ are disjoint, we have the estimate

$$\|T(\chi_{L_1^c} f)\|_{L^q(K_0)} \leq \|A\|_\infty\, C_{K_0, L_1} \|f\|_{L^2(\mathbb{R}^N)} \tag{2.2}$$

Indeed, let $\tilde{K}$ be so that $K_0 \subset\subset \tilde{K} \subset\subset L_1$. Then by the boundedness of the Riesz transform on $L^q(\mathbb{R}^N)$,

$$\begin{aligned}
\|T(\chi_{L_1^c} f)\|_{L^q(K_0)} \leq & \|R_i(\chi_{\tilde{K}} A_{ij} R_j((\chi_{L_1^c} f))\|_{L^q(K_0)} + \|R_i(\chi_{\tilde{K}^c} A_{ij} R_j((\chi_{L_1^c} f))\|_{L^q(K_0)}\\
\leq & \|A\|_{L^\infty(\mathbb{R}^N)} \|R_j(\chi_{L_1^c} f)\|_{L^q(\tilde{K})} + \|R_i(\chi_{\tilde{K}^c} A_{ij} R_j((\chi_{L_1^c} f))\|_{L^q(K_0)}
\end{aligned}$$

We now apply the Cauchy-Schwarz inequality to obtain

$$\begin{aligned}
\|R_j(\chi_{L_1^c} f)\|_{L^q(\tilde{K})} = & \left( \int_{K_0} \left| \int_{\mathbb{R}^N \setminus L_1} f(y) \frac{x_j - y_j}{|x-y|^{N+1}}\, dy \right|^q dx \right)^{\frac{1}{q}}\\
\leq & \left( \int_{K_0} \|f\|_{L^2(\mathbb{R}^N)}^q \left( \int_{\mathbb{R}^N \setminus L_1} \frac{1}{|x-y|^{2N}}\, dy \right)^{q/2} dx \right)^{\frac{1}{q}}\\
\leq & C|K_0|^{1/q} \|f\|_{L^2(\mathbb{R}^N)} \left( \int_c^\infty \frac{1}{t^{2N}} t^{N-1}\, dt \right)^{\frac{1}{2}}\\
\leq & C_{K_0, L_1, q} \|f\|_{L^2(\mathbb{R}^N)},
\end{aligned}$$

where we have used the disjointness of $K_0$ and $L_1^c$ (in particular that $\operatorname{dist}(K_0, L_1^c) = c > 0$). A similar argument shows that

$$\|R_i(\chi_{\tilde{K}^c} A_{ij} R_j((\chi_{L_1^c} f))\|_{L^q(K_0)} \leq C_{\tilde{K}, L_1, q} \|A_{ij} R_j((\chi_{L_1^c} f))\|_{L^2(\mathbb{R}^N)},$$

and so using the boundedness of the Riesz transform on $L^2(\mathbb{R}^N)$, we conclude that

$$\|R_i(\chi_{\tilde{K}^c} A_{ij} R_j((\chi_{L_1^c} f))\|_{L^q(K_0)} \leq C_{\tilde{K}, K_0, q} \|A\|_{L^\infty(\mathbb{R}^N)} \|f\|_{L^2(\mathbb{R}^N)}.$$

It thus remains to estimate $II$. Let us begin by observing that the commutator estimates with a Lipschitz continuous function (see Proposition 2.3) imply that

$$\begin{aligned}
II = & \|\mathcal{C}(\eta, T)[\chi_{L_1} f]\|_{L^q(\mathbb{R}^N)}\\
\leq & C_\eta (\| I_1 |\chi_{L_1} f| \|_{L^q(\mathbb{R}^N)} + \| I_1 |A_{ij} R_j(\chi_{L_1} f)| \|_{L^q(\mathbb{R}^N)}).
\end{aligned}$$

In particular, $q > 2$ implies that $Nq/(N+q) > 1$ and so $I_1 : L^{Nq/(N+q)}(\mathbb{R}^N) \to L^q(\mathbb{R}^N)$ is bounded. Moreover, $R_j : L^r(\mathbb{R}^N) \to L^r(\mathbb{R}^N)$ is bounded for $1 < r < +\infty$, which combined with the fact that $A \in L^\infty(\mathbb{R}^N; \mathbb{R}^{N \times N})$ (recall that $N \geq 2$) implies that

$$II \leq C \|f\|_{L^{Nq/(N+q)}(L_1)}.$$

If we let $L_0 := \Omega_1$, then our estimates show that

$$\|f\|_{L^{q_0}(L_0)} \leq C \left( \|T(f)\|_{L^{q_0}(K_0)} + \|f\|_{L^2(\mathbb{R}^N)} + \|f\|_{L^{q_1}(L_1)} \right)$$

for $q_i := Nq/(N + iq)$. Now, if $q_1 \leq 2$ then an application of Hölder's inequality implies the desired result. Otherwise we iterate the previous argument by finding

$$K_0 \subset\subset L_1 \subset\subset K_1 \subset\subset L_2 \subset\subset \ldots K_i \subset\subset L_{i+1} \subset\subset \Omega_2$$

to obtain the estimate

$$\|f\|_{L^{q_i}(L_i)} \leq C \left( \|T(f)\|_{L^{q_i}(K_i)} + \|f\|_{L^2(\mathbb{R}^N)} + \|f\|_{L^{q_{i+1}}(L_{i+1})} \right),$$

provided $q_{i+1} > 1$ (in order that $I_1 : L^{q_{i+1}}(\mathbb{R}^N) \to L^{q_i}(\mathbb{R}^N)$). However, $q_i > 2$ implies $q_{i+1} > 1$, and so we continue the iteration a finite number of times until we obtain that $q_j \leq 2$ for some $j \in \mathbb{N}$. Then collecting the terms our estimate reads

$$\|f\|_{L^q(\Omega_1)} \leq C \left( \sum_{i=0}^{j-1} \|T(f)\|_{L^{q_i}(K_i)} + \|f\|_{L^2(\mathbb{R}^N)} + \|f\|_{L^{q_j}(L_j)} \right),$$

from which the inequality (2.2) is a simple consequence of Hölder's inequality, and thus the proposition is established. ∎

Finally, we require the following result.

**Proposition 2.4** *Let $\Omega \subset \mathbb{R}^N$ be open and bounded, $s \in [0, N)$, and $2 \leq p < +\infty$. Assume that for all $\varphi \in C_c^\infty(\Omega)$,*

$$\int f(-\Delta)^{\frac{s}{2}} \varphi = \int h(-\Delta)^{\frac{s}{2}} \varphi.$$

*Then for $\Omega_1 \subset\subset \Omega$, there exists a constant $C = C(\Omega_1)$ such that*

$$\|f\|_{L^p(\Omega_1)} \leq C \left( \|h\|_{L^p(\mathbb{R}^N)} + \|f\|_{L^2(\mathbb{R}^N)} \right).$$

**Proof.** Let $\Omega_1 \subset\subset \Omega_2 \subset\subset \Omega$ and $\varphi \in C_c^\infty(\Omega_2)$ be such that

$$\|f\|_{L^p(\Omega_1)} \leq 2 \int f \varphi$$

and $\|\varphi\|_{L^{p'}(\mathbb{R}^N)} \leq 1$.

We argue by first reducing to the case where the support of $\varphi$ is a ball. We can accomplish this by covering $\Omega_2$ with finitely many balls $B(x_j, r_j)$ of controlled overlap such that $B(x_j, 4r_j) \subset\subset \Omega$, where the number of balls can be taken to depends only on the distance of $\Omega_1$ to $\Omega^c$. Then by subordinating a partition of unity to balls $B(x_j, r_j)$ we can write

$$\varphi = \sum_{j=1}^{l} \varphi_j$$

with $\mathrm{supp}\,\varphi_j \subset B(x_j, r_j)$ for each $j$ and $|\varphi_j| \leq |\varphi|$. Then for $j$ fixed we have

$$\int f\, \varphi_j = 2\int f(-\Delta)^{\frac{s}{2}} I_s \varphi_j$$

$$= 2\int f(-\Delta)^{\frac{s}{2}}(\eta_j I_s \varphi) + 2\int f(-\Delta)^{\frac{s}{2}}((1-\eta_j)I_s \varphi_j)$$

$$= 2\int h(-\Delta)^{\frac{s}{2}}(\eta_j I_s \varphi_j) + 2\int f(-\Delta)^{\frac{s}{2}}((1-\eta_j)I_s \varphi_j)$$

$$\leq 2\left(\|h\|_{L^p(\mathbb{R}^N)}\|(-\Delta)^{\frac{s}{2}}(\eta_j I_s \varphi)\|_{L^{p'}(\mathbb{R}^N)} + \|f\|_{L^2(\mathbb{R}^N)}\|(-\Delta)^{\frac{s}{2}}((1-\eta_j)I_s \varphi)\|_{L^2(\mathbb{R}^N)}\right),$$

where $\eta_j \in C_c^\infty(\Omega)$ with $\eta \equiv 1$ on $B(x_j, 4r_j)$. Then if we can establish the estimates

$$\|(-\Delta)^{\frac{s}{2}}(\eta_j I_s \varphi)\|_{L^{p'}(\mathbb{R}^N)} \leq C\|\varphi_j\|_{L^{p'}(\mathbb{R}^N)} \tag{2.3}$$

$$\|(-\Delta)^{\frac{s}{2}}((1-\eta_j)I_s \varphi)\|_{L^2(\mathbb{R}^N)} \leq C\|\varphi_j\|_{L^{p'}(\mathbb{R}^N)}, \tag{2.4}$$

the result will follow by summing in $j$ and using the pointwise inequality $|\varphi_j| \leq |\varphi|$.

Let us therefore first examine (2.3), and to save notation we drop the dependence in $j$. If we take the three term commutator $H_s$ introduced by Da Lio and Rivière [3]

$$H_s(\eta, I_s\varphi) := (-\Delta)^{\frac{s}{2}}(\eta I_s \varphi) - (-\Delta)^{\frac{s}{2}}\eta\, I_s\varphi - \eta\varphi,$$

we can use

$$\|H_s(\eta, I_s\varphi)\|_{L^{p'}(\mathbb{R}^N)} \leq C\|\varphi\|_{L^{p'}(\mathbb{R}^N)}.$$

This estimate follows via the Hardy-Littlewood decomposition in [3] or using the *pointwise* estimates in [10] (see [4, Theorem 1.2] for a precise version that can be applied here and also [1, 2] for various extensions). Thus, it suffices to show that

$$\|(-\Delta)^{\frac{s}{2}}\eta\, I_s\varphi\|_{L^{p'}(\mathbb{R}^N)} + \|\eta\varphi\|_{L^{p'}(\mathbb{R}^N)} \leq C\|\varphi\|_{L^{p'}(\mathbb{R}^N)}.$$

The second term can be estimated in terms of the right hand side trivially since $|\eta| \leq 1$, while for the first term one applies Hölder's inequality with exponent $Np'/(N-sp')$ and its Hölder conjugate $r$ when $N - sp' > 0$ (Note that from $\eta \in C_c^\infty(\mathbb{R}^n)$ we know that $(-\Delta)^{\frac{s}{2}}\eta \in L^r(\mathbb{R}^n)$ for any $r \in (1, \infty)$ e.g. by interpolation.), which yields

$$\|(-\Delta)^{\frac{s}{2}}\eta I_s\varphi\|_{L^{p'}(\mathbb{R}^N)} \leq \|(-\Delta)^{\frac{s}{2}}\eta\|_{L^r(\mathbb{R}^N)}\|I_s\varphi\|_{L^{Np'/(N-sp')}(\mathbb{R}^N)}$$

$$\leq C\|\varphi\|_{L^{p'}(\mathbb{R}^N)}.$$

If $N - sp' < 0$, then

$$\|(-\Delta)^{\frac{s}{2}}\eta I_s\varphi\|_{L^{p'}(\mathbb{R}^N)} \leq \|(-\Delta)^{\frac{s}{2}}\eta\|_{L^{p'}(\mathbb{R}^N)}\|I_s\varphi\|_{L^\infty(\mathbb{R}^N)}$$
$$\leq C\|\varphi\|_{L^{p'}(\mathbb{R}^N)}$$

follows from the fact that $\varphi$ has compact support. When $N - sp' = 0$, we take $\tilde{p}' < p'$ and set $\frac{1}{\tilde{r}} := \frac{1}{p'} - \frac{1}{\tilde{p}'}$, then

$$\|(-\Delta)^{\frac{s}{2}}\eta I_s\varphi\|_{L^{p'}(\mathbb{R}^N)} \leq \|(-\Delta)^{\frac{s}{2}}\eta\|_{L^{\tilde{r}}(\mathbb{R}^N)}\|I_s\varphi\|_{L^{N\tilde{p}'/(N-s\tilde{p}')}(\mathbb{R}^N)}$$
$$\leq C\|\varphi\|_{L^{\tilde{p}'}(\mathbb{R}^N)}.$$

The estimate follows again in this case by the fact that $\varphi$ has compact support.

Finally, to establish (2.4) we write

$$(1 - \eta) = \sum_{k=2}^{\infty} \theta_{A_{2^k r}},$$

where each $\theta_{A_{2^k r}}$ is supported on an annulus of width $2^k r$. Then disjoint support arguments (see, for example, Lemma 3.7 in [9]) imply the estimate

$$\|(-\Delta)^{\frac{s}{2}}(\theta_{A_{2^k r}} I_s\varphi)\|_{L^2(\mathbb{R}^N)} \leq C(2^k r)^{-N/2} r^{N/p}\|\varphi\|_{L^{p'}(\mathbb{R}^N)},$$

from which we obtain

$$\|(-\Delta)^{\frac{s}{2}}((1-\eta) I_s\varphi)\|_{L^2(\mathbb{R}^N)} \leq \sum_{k=2}^{\infty} \|(-\Delta)^{\frac{s}{2}}(\theta_{A_{2^k r}} I_s\varphi)\|_{L^2(\mathbb{R}^N)}$$
$$\leq \left( C \sum_{k=2}^{\infty} (2^k r)^{-N/2} r^{N/p} \right) \|\varphi\|_{L^{p'}(\mathbb{R}^N)}.$$

As the series is summable we have established the desired inequality and therefore the theorem is proved. ∎

We are now ready to prove the main result.

**Proof of Theorem 1.2.** Suppose $G \in L^p(\mathbb{R}^N; \mathbb{R}^N)$ and $u \in H^s(\mathbb{R}^N)$ satisfies the equation (1.4). The claim of this theorem is that for any $K \subset\subset \Omega$, one has the estimate

$$\|D^s u\|_{L^p(K;\mathbb{R}^N)} \leq C\left(\|G\|_{L^p(\mathbb{R}^N;\mathbb{R}^N)} + \|(-\Delta)^{\frac{s}{2}} u\|_{L^2(\mathbb{R}^N)}\right).$$

We will see that the result is a consequence of a combination of Propositions 2.2 and 2.4, and we argue as follows. Define $g := R^* G = -\sum_{j=1}^N R_j G_j$, so that $g \in L^p(\mathbb{R}^N)$ and $u$ satisfies

$$\int_\Omega T(-\Delta)^{\frac{s}{2}} u \, (-\Delta)^{\frac{s}{2}} \varphi = \int g \, (-\Delta)^{\frac{s}{2}} \varphi \quad \forall \varphi \in C_c^\infty(\Omega),$$

where $T$ is as in Proposition 2.1. Moreover, a cutoff argument similar to those previously employed implies that if $K \subset\subset \Omega_1$, then one has

$$\|D^s u\|_{L^p(K;\mathbb{R}^N)} = \|R(-\Delta)^{\frac{s}{2}} u\|_{L^p(K;\mathbb{R}^N)}$$
$$\leq C\left(\|(-\Delta)^{\frac{s}{2}} u\|_{L^p(\Omega_1)} + \|(-\Delta)^{\frac{s}{2}} u\|_{L^2(\mathbb{R}^N)}\right),$$

and so this and boundedness of the Riesz transforms (to obtain bounds on $g$ in terms of $G$ in $L^p$) imply that it suffices to show the estimate

$$\|(-\Delta)^{\frac{s}{2}}u\|_{L^p(\Omega_1)} \leq C \left(\|g\|_{L^p(\mathbb{R}^N)} + \|(-\Delta)^{\frac{s}{2}}u\|_{L^2(\mathbb{R}^N)}\right).$$

for $\Omega_1 \subset\subset \Omega$.

We first apply Proposition 2.2 with $f = (-\Delta)^{\frac{s}{2}}u$ and for $\Omega_1 \subset\subset \Omega_2 \subset\subset \Omega$ yielding

$$\|(-\Delta)^{\frac{s}{2}}u\|_{L^p(\Omega_1)} \leq C \left(\|T(-\Delta)^{\frac{s}{2}}u\|_{L^p(\Omega_2)} + \|(-\Delta)^{\frac{s}{2}}u\|_{L^2(\mathbb{R}^N)}\right).$$

Now Proposition 2.4 and boundedness of $T: L^2(\mathbb{R}^N) \to L^2(\mathbb{R}^N)$ gives

$$\|T(-\Delta)^{\frac{s}{2}}u\|_{L^p(\Omega_2)} \leq C \left(\|g\|_{L^p(\mathbb{R}^N)} + \|T(-\Delta)^{\frac{s}{2}}u\|_{L^2(\mathbb{R}^N)}\right)$$
$$\leq C \left(\|g\|_{L^p(\mathbb{R}^N)} + \|(-\Delta)^{\frac{s}{2}}u\|_{L^2(\mathbb{R}^N)}\right).$$

Therefore, we find

$$\|(-\Delta)^{\frac{s}{2}}u\|_{L^p(\Omega_1)} \leq C \left(\|g\|_{L^p(\mathbb{R}^N)} + \|(-\Delta)^{\frac{s}{2}}u\|_{L^2(\mathbb{R}^N)}\right),$$

which is the thesis. ∎

## 3 Acknowledgements

The first author is supported by the SNF and would like to thank D. Spector, NCTU and the CMMSC for their hospitality during the work of this project. The second author is partially supported by the Taiwan Ministry of Science and Technology under research grant MOST 103-2811-M-008-066 and would like to thank Prof. Chern Jann-Long and Prof. Chern I-Liang for their support. The third author is supported by the Taiwan Ministry of Science and Technology under research grant MOST 103-2115-M-009-016-MY2.